\documentclass{amsart}

\usepackage{amsmath, amssymb, amsthm, amsfonts, marginnote, stmaryrd}

\usepackage{comment}

\usepackage{xcolor}

\input xy
\xyoption{all}

\usepackage{hyperref}

\swapnumbers
\numberwithin{equation}{section}

\theoremstyle{plain}
\newtheorem{theorem}[subsubsection]{Theorem}

\newtheorem{prop}[subsubsection]{Proposition}

\newtheorem{conj}[subsubsection]{Conjecture}

\newtheorem{ansatz}[subsubsection]{Ansatz}

\theoremstyle{definition}
\newtheorem{defn}[subsubsection]{Definition}

\newtheorem{remark}[subsubsection]{Remark}

\newtheorem{ex}[subsubsection]{Example}

\def\enhIndCoh{\textup{IndCoh}^{\diamondsuit}}

\def\module{\operatorname{-mod}}

\def\cat{\operatorname{-cat}}

\def\module{\mathrm{mod}}

\def\AA{\mathbb{A}}

\def\CC{\mathbb{C}}

\def\FF{\mathbb{F}}
\def\GG{\mathbb{G}}

\def\HH{\mathbb{H}}

\def\LL{\mathbb{L}}

\def\PP{\mathbb{P}}

\def\SS{\mathbb{S}}

\def\ZZ{\mathbb{Z}}
\def\BZ{\mathbb{Z}}

\newcommand\cA{\mathcal{A}}
\newcommand\cB{\mathcal{B}}
\newcommand\cC{\mathcal{C}}
\newcommand\cD{\mathcal{D}}

\newcommand\cH{\mathcal{H}}

\newcommand\cL{\mathcal{L}}
\newcommand\cM{\mathcal{M}}
\newcommand\cN{\mathcal{N}}
\newcommand\cO{\mathcal{O}}

\newcommand\cU{\mathcal{U}}

\newcommand\cZ{\mathcal{Z}}

\newcommand\frg{\mathfrak{g}}
\newcommand\fg{\mathfrak{g}}

\newcommand\frt{\mathfrak{t}}

\newcommand\frtv{\check{\frt}}

\newcommand{\Conn}{\textup{Conn}}

\newcommand{\Ind}{\textup{Ind}}

\newcommand\Loc{\textup{Loc}}

\newcommand\Mod{\textup{Mod}}

\newcommand{\Perf}{\textup{Perf}}

\newcommand{\Pic}{\textup{Pic}}

\newcommand{\QCoh}{\textup{QCoh}}

\newcommand{\reg}{\textup{reg}}

\newcommand\Spec{\textup{Spec}}

\newcommand\St{\mathit{St}}

\newcommand\Sym{\textup{Sym}}

\newcommand{\Tr}{\textup{Tr}}
\newcommand{\tr}{\textup{tr}}

\newcommand\Hom{\textup{Hom}}
\newcommand\End{\textup{End}}
\newcommand\Map{\textup{Map}}

\newcommand{\Gm}{\GG_m}

\renewcommand\a\alpha
\renewcommand\b\beta
\newcommand\g\gamma
\newcommand\G\Gamma
\renewcommand\d\delta
\newcommand\D\Delta

\newcommand{\wh}[1]{\widehat{#1}}
\newcommand\quash[1]{}

\newcommand\bs{\backslash}
\newcommand\ot{\otimes}

\newcommand{\Coh}{\textup{Coh}}

\newcommand{\beq}{\begin{equation}}
\newcommand{\eeq}{\end{equation}}

\newcommand{\Ch}{\mathit{Ch}}

\newcommand{\gitquot}{\hspace{-0.1em}/\hspace{-0.25em}/}

\newcommand{\Gv}{{\check{G}}}

\newcommand{\Tv}{{\check{T}}}
\newcommand{\Xv}{\check{X}}
\newcommand{\Bv}{\check{B}}

\newcommand{\Mv}{\check{M}}
\newcommand{\sv}{\check{s}}

\newcommand{\Cx}{{\CC^\times}}

\newcommand{\actson}{\circlearrowright}

\usepackage[mathscr]{euscript}

\newcommand{\IndCoh}{\textup{IndCoh}}
\newcommand{\Tors}{\textup{Tors}}

\newcommand{\Tate}{\textup{Tate}}

\newcommand{\GIT}{{/\! /}}

\newcommand{\AlgStk}{\textup{AlgStk}}

\title[Potent  Categorical Representations]{Potent Categorical Representations}

\author{David Ben-Zvi} \address{Department of Mathematics\\University
  of Texas\\Austin, TX 78712-0257} \email{benzvi@math.utexas.edu}
\author{David Nadler} \address{Department of Mathematics\\University
  of California\\Berkeley, CA 94720-3840}
\email{nadler@math.berkeley.edu}

\dedicatory{Dedicated to Tom Nevins}
\thanks{Supported by NSF individual grants DMS-2001398, DMS-2302356 (DBZ) and DMS-2101466, DMS-2401178 (DN)}

\begin{document}


\begin{abstract}
We introduce and motivate -- based on ongoing joint work with Germ\'an Stefanich -- the notion of {\em potent categorical representations} of a complex reductive group $G$, specifically a conjectural Langlands correspondence identifying potent categorical representations of $G$ and its Langlands dual $\Gv$.
We emphasize the symplectic nature of potent categorical representations in their simultaneous dependence on parameters in  maximal tori for $G$ and $\Gv$, specifically how their conjectural Langlands correspondence fits within a 2-categorical  Fourier  transform. Our key tool to make  various ideas precise is higher sheaf theory and its microlocalization, specifically a theory of ind-coherent sheaves of categories on stacks. The constructions are inspired by the physics of 3d mirror symmetry and S-duality on the one hand, and the theory of double affine Hecke algebras on the other.  We also highlight further conjectures 
related to ongoing programs in and around geometric representation theory. 
\end{abstract}

\maketitle


\section{Introduction}

In this  brief introduction, we indicate what ``categorical representations"  we will encounter and  why we use the term ``potent" for the developments we propose. The rest of the paper is split into three sections aimed at different audiences. Section \S\ref{GRT intro} introduces the main ideas from the point of view of geometric representation theory and derived algebraic geometry. Section \S\ref{TFT section} recasts the prior constructions within the framework of topological quantum field theory. 
Section \S\ref{applications} emphasizes links to the representation theory of Hecke algebras and presents some conjectural applications to a variety of topics including the local Langlands correspondence, Hilbert schemes and homology of character varieties. All of the new material presented here has been developed with Germ\'an Stefanich and will be further detailed in our forthcoming joint work~\cite{BZNS}.


\subsection{What are ``categorical representations''?} 

A core topic in representation theory is the study of modules
for group algebras, such as distributions on finite, Lie or $p$-adic groups $G$. These modules naturally arise in the form of generalized functions or cohomology of $G$-spaces $X$. 
One typical meaning of {\em categorical representation theory} 
is the parallel study of module categories for categorical group algebras -- monoidal categories of sheaves of one flavor or another on groups $G$. Such categorical representations naturally arise as categories of sheaves of the same flavor on $G$-spaces $X$.

In both of the above pursuits, a key choice is what kinds of functions or sheaves are considered. The developments of this article  will most closely encounter {\em de Rham group algebras} of reductive groups $G$ where one takes the
 categorical group algebra to be the dg category of $\cD$-modules $\cD(G)$, and fundamental examples of representations are the dg categories $\cD(G/K)$  of  $\cD$-modules on homogenous spaces and $U\fg-\module$ of representations of the Lie algebra $\fg$ (see  \cite{dhilloninformal} for an excellent survey).
So here we are in the realm of Hecke categories (as often studied via Soergel bimodules, see~\cite{soergelbimodules}) which
form the intertwiners on induced representations, 
 and the local geometric Langlands program, which classifies de Rham categorical representation of the loop group~$LG$.

\subsection{What is ``potent''?} We use the term {\em potent} to contrast with {\em unipotent}, just as {\em crepant} is used  to contrast with {\em discrepant}. In traditional representation theory of reductive groups $G$, induced representations come in families (or series) parametrized by Weyl group orbits of characters $\chi$ of various tori $T$. The term unipotent
refers to  those representations whose parameter $\chi$ is the same as that of the trivial representation. 

In de Rham categorical representation theory, 
induced representations likewise come in families parametrized by Weyl group orbits, but now of local systems $\cL_\chi$ on tori $T$. 
We use the term {\em potent} to signify that we will not only allow these parameters $\cL_\chi$ to vary, but take seriously the geometry of their parameter space $\check T$. We will see that this naturally leads to the {\em microlocal geometry} of the parameter space~$\check T$, so that in fact there will be dual momentum parameters associated with $T$. 
Alternatively, or in fact equivalently under duality,  the parameters associated with $T$ will also arise as equivariant parameters to keep track of {\em genuine equivariant topology}.
In the most symmetric realization, potent categorical representations will  live over the symplectic space $T \times \check T$ and be   parametrized by  Lagrangians therein.


\subsection{Conventions}\label{conventions}
In what follows, all algebra is understood to be over the complex numbers $\CC$ or at most a field $k$ of characteristic 0. All category theory is understood to be derived in an $\infty$-categorical sense, and we often implicitly work within the $\infty$-cateogory $\Pr^L_k$ of presentable $k$-linear $\infty$-categories or $\St_k$ of  presentable stable $k$-linear 
$\infty$-categories.
 
 We have tried to balance intuition and precision, but there are specific aspects where care in interpretation is particularly warranted.
 First, many of our constructions have alternative  realizations that differ algebraically but are the same analytically. For example, many constructions have Betti and de Rham versions associated  to either an algebraic torus $T$ or the quotient stack $\frt/\Lambda_T$ of its Lie algebra by the coweight lattice. We have tried to cleanly convey broad themes without dwelling on this difference, while then stating specific results in as precise a form as possible. We hope this provides more clarity than confusion. 
Second, many of our constructions lead to two-periodic categories, or more generally, categories over some additional parameters.  We are not particularly careful to make this explicit and leave implicit requisite base-changes of other categories to these additional parameters.





\subsection{Acknowledgements}
This article is based on ongoing joint work with Germ\'an Stefanich, which in turn is based on his work on higher sheaf theory. We are grateful to him for letting us share key ideas here and  providing valuable feedback on the exposition. Of course, any errors are due to the named authors. 

We are  grateful  to the following collaborators, colleagues, and  researchers who have generously shared their ideas and understanding with us:  Andrew Blumberg, Sanath Devalapurkar, Davide Gaiotto, Ben Gammage, Rok Gregoric, Sam Gunningham, Justin Hilburn, Matt Hogancamp, Quoc Ho, David Jordan, Penghui Li, Andy Neitzke, Toly Preygel, Pavel Safronov, Jeremy Taylor, Akshay Venkatesh, Edward Witten, Enoch Yiu.


\section{Perspective of geometric representation theory}\label{GRT intro}

This section is an 
overview of potent categorical representations from the perspective of geometric representation theory. 
We start with a review of $\cD$-modules via loops spaces as motivation, then give an informal definition of potent categorical representations.
We then sketch what is needed to make the informal definition 
more precise,
and finally  explore its calculation and duality via Hecke theory.


\subsection{Loop spaces and $\cD$-modules}\label{loops and conns section}

The {\em derived loop space} of a scheme or stack $X$ is the mapping space
$\cL X=\Map(S^1,X)=X\times_{X\times X} X,$ where $S^1=B\ZZ$ is the simplicial (or ``animated") circle. 

For $X$ a scheme with cotangent complex $\LL_X$ and shifted tangent bundle $T[-1]X \simeq \Spec_X \Sym_X(\LL_X[1])$, the HKR theorem 
provides an exponential map identifying $T[-1]X$ with $\cL X$.
The $S^1$-action on $\cL X$ by loop rotation is encoded by the de Rham differential $d$ viewed as a degree $-1$ vector field. Passing to cyclic functions gives the identification of periodic cyclic homology with de Rham cohomology due to Connes and Feigin-Tsygan.

\subsubsection{$\cD$-modules on schemes} For $X$ a scheme, 
it was shown in~\cite{conns, TV1, TVChern} (for smooth schemes) and ~\cite{preygelloops} (in general) that two-periodic $\cD$-modules on $X$ are  equivalent to cyclic sheaves on the loop space $\cL X$:
\begin{equation}\label{eq:loops D}
    \cD(X)\ot k[u,u^{-1}]\simeq \IndCoh(\cL X)^{\Tate}
\end{equation}

Here for $Z$  a quasi-smooth scheme, $\IndCoh(Z)$ (a categorified analog of algebraic distributions) assigns the category of ind-coherent sheaves, an enlargement of the category of quasicoherent sheaves $\QCoh(Z)$ (a categorified analog of algebraic functions) that allows more refined phenomena at singularities (non-trivial wavefront sets). More precisely, 
following~\cite{ArinkinGaitsgory}, objects $\cM\in \IndCoh(Z)$ have a singular support inside $T^*[1] Z$, 
and the image of the inclusion  $\QCoh(Z) \subset \IndCoh(Z)$ is
objects with singular support in the zero-section. 
Furthermore, there is the precise relation under traces (i.e.~Hochschild homology) $\tr(\IndCoh(Z))\simeq \omega(\cL Z)$ and $\tr(\QCoh(Z))\simeq \cO(\cL Z)$. 

The superscript ``$\Tate$" denotes that we perform a Tate (i.e.~periodic cyclic) construction: we first pass to the (homotopy) $S^1$-invariant category $ \IndCoh(\cL X)^{S^1}$, which is linear over the equivariant cohomology $H^\ast(BS^1)\simeq k[u]$, $|u| = 2$, 
and then invert the equivariant parameter $u$ to obtain $\IndCoh(\cL X)^{\Tate} =
 \IndCoh(\cL X)^{S^1}[u^{-1}]$. This will be the first of many times we  implement a Tate construction, and we flag here a key theme from equivariant homotopy theory encountered throughout: 
 to be sure  a Tate construction is non-trivial, one needs to take care the notion of $S^1$-action is rich enough to see more than $u$-torsion.
 
 The equivalence~\eqref{eq:loops D} of the Tate construction with two-periodic $\cD$-modules goes via Koszul duality, under which this $u$-deformation of $\cL X \simeq T[-1] X$ corresponds to the $\hbar$-deformation quantization of $T^*X$. There are many further compatibilities of the equivalence~\eqref{eq:loops D}, for example between the singular supports of $\cD$-modules and that of ind-coherent sheaves.

\subsubsection{$\cD$-modules on stacks}\label{potent D-module section} 
For $X$ a stack, the loop space $\cL X$ is a rich global object which can not be recovered from the loop spaces of a cover $U \to X$. 

\begin{ex}\label{ex:glob quot}
    For a global quotient $X=Z/G$, the loop space is the {\em derived inertia stack} 
$\cL (Z/G) \simeq \{z\in Z, g\in G_z\}/G
$
where $G_z \subset G$ is the stabilizer of $z\in Z$.
When $G$ is reductive, the natural map $\cL (Z/G) \to \cL BG = G/G \to T\gitquot W$ allows for a
 {\em Jordan decomposition of loops}~\cite{reps,Harrison,BZCHN2}.
\end{ex}

So we arrive at a divergence: $\cD$-modules on a stack $X$ are defined so as to satisfy descent for a cover, but cyclic sheaves on the loop space $\cL X$ do not.
On the one hand, it is possible to identify $\cD$-modules within all cyclic sheaves: 
for $X$ a smooth stack, 
it was shown in~\cite{conns}  that  two-periodic $\cD$-modules are  equivalent to cyclic sheaves on the formal completion  $\wh \cL X \subset \cL X$ along constant loops $X\subset \cL X$:
$$\cD(X)\ot k[u,u^{-1}]\simeq \IndCoh(\wh \cL X)^{\Tate}$$

 On the other hand, cyclic sheaves on the entire loop space provide  a richer object to study  which incorporates equivariant parameters (such as the familiar base $T\gitquot W$ for $G$-equivariant $K$-theory as seen in Example~\ref{ex:glob quot}). Since they play an important role in our later developments, we will distinguish them with a name consistent with those developments. 

\begin{defn} {\em Potent $\cD$-modules} on a stack $X$ are cyclic sheaves on its loop space
$$\cD^{pot}(X):=\IndCoh(\cL X)^{\Tate}$$
\end{defn}
 
Potent $\cD$-modules are closely related to reduced $K$-motives~\cite{eberhardtK,eberhardteteve,tayloruniversal}, which map to potent $\cD$-modules via the Chern character~\cite{TVChern,HSS}.


For  quotients of flag varieties $X=K\bs G/B$ by subgroups $K \subset G$, we argued in~\cite{reps, BZCHN, BZCHN2} that derived loop spaces and potent $\cD$-modules form a natural setting for the local Langlands correspondence in both archimedean and non-archimedean settings.


\subsection{Potent categorical representations}\label{sect:pcrs}
We will implement the analogy: $\cD$-modules are to potent $\cD$-modules as de Rham categorical representations are to potent categorical representations. In this analogy, the notion of sheaves needs to be replaced with a suitable theory of sheaves of categories. Our main tool
 will be higher sheaf theory as developed by Stefanich~\cite{stefanichPr,stefanichIndCoh, Scategorification}. 

\subsubsection{$2\IndCoh$}\label{sect:2indcoh}

The theory of quasicoherent sheaves of categories
$2\QCoh$ as developed in~\cite{1affine}  gives a robust  
categorified analog of
quasicoherent sheaf theory (itself a  categorified analog of algebraic functions). Stefanich's work provides a notion of ``ind-coherent sheaves of categories'' $2\IndCoh$ which is a categorified analog of ind-coherent sheaf theory $\IndCoh$ (itself a  categorified analog of of algebraic distributions). 
Among many structures, 
the theory assigns a 2-category $2\IndCoh(X)$ to a  stack $X$
(typically assumed to be smooth)  in a functorial way so that maps of stacks $Y \to X $ (typically assumed to be proper) provide objects $\cM(Y)$ and morphisms $\cM(Y_1) \to \cM(Y_2)$ are calculated by $\IndCoh(Y_1 \times_X Y_2)$. 
Objects $\cM \in 2\IndCoh(X) $ have a singular support inside $T^*[2] X$, and in analogy with the  inclusion $\QCoh \subset \IndCoh$, there is an inclusion  $2\QCoh \subset 2\IndCoh$ whose image is objects with singular support in the zero-section. Furthermore, there is the precise relation under traces $\Tr(2\IndCoh(X))\simeq \IndCoh(\cL X)$ and $\Tr(2\QCoh(X))\simeq \QCoh(\cL X)$. So for $X$ a scheme, the cyclic trace recovers two-periodic $\cD$-modules 
$\Tr(2\IndCoh(X))^{\Tate}\simeq \IndCoh(\cL X)^{\Tate}\simeq \cD(X) \otimes k[u,u^{-1}]$. 



 \begin{remark}\label{rem:2indcoh has many generators}
 We will refer to objects 
of $2\IndCoh(X)$ as {\em ind-coherent sheaves of categories} but
make the following warning even  in the case when $X$ is affine: On the one hand, objects of $2\QCoh(X)$ are sheaves of categories  over $X$ in the Zariski topology and in fact determined by their global sections $2\QCoh(X)=\QCoh(X)-\module$. On the other hand, objects of $2 \IndCoh(X)$ encode far more data: to view them as sheaves of categories requires far more test schemes $Y \to X$, and in particular, they are not determined by their global sections
$\Gamma(X, Y) = \Hom(X,Y) = \IndCoh(Y)$
in any non-trivial situation.
For example, there are non-trivial objects with vanishing global sections. This divergence between $2\QCoh$ and $2\IndCoh$, or in fact $\QCoh$ and $\IndCoh$,  reflects  the phenomenon of {\em expansion} (familiar in equivariant stable homotopy theory under the name {\em decomposition}) discussed below in \S\ref{expansion section}.
It leads to the phenomenon that  the 2-categories appearing in this article  are typically not realized as modules for a monoidal category. 

 \end{remark}

 We will need  two  key elaborations on $2\IndCoh$ --  {\em periodization} and {\em equivariant expansion}, as discussed below in \S\ref{IndCoh elaborations} 
 -- whose outcome we will denote by $2\enhIndCoh$. For the moment, we will proceed here with $2\enhIndCoh$ and the reader can either blackbox it as an elaboration on $2\IndCoh$ or first navigate  to \S\ref{IndCoh elaborations} and then return here.


\subsubsection{Potent categorical representations}\label{def of potent reps section}
Now let's return to 
$G$  a reductive group, and $\cD(G)$ the de Rham group algebra of $\cD$-modules on $G$. We regard  $\cD(G)$ as an algebra object in presentable stable $k$-linear  $\infty$-categories.
The 2-category of de Rham categorical $G$-representations is  given by modules for the de Rham group algebra 
$G-cat := \cD(G)-\module.
$
We would like to trade the algebra of $G$ for the geometry of its classifying stack $BG$. 
Via 
the theory of 
the 
prior section, we can reinterpret de Rham categorical $G$-representations via loops in the following form
$
G-{cat}  \simeq 2\IndCoh(\wh G/G)^{\Tate}.
$
Here $\wh \cL BG = \wh G/G \subset G/G = \cL BG$ is the formal loop space, and we take cyclic objects with respect to loop rotation.

In analogy with potent $\cD$-modules, we arrive at
potent categorical $G$-representations by extending to cyclic objects on the entire  loop space:

\begin{defn}
{\em  Potent categorical $G$-representations} are cyclic ind-coherent sheaves of categories 
on the loop space 
$$
G-{cat}^{pot} := 2\enhIndCoh(G/G)^{\Tate}
$$
\end{defn}

Just as potent $\cD$-modules on $X$ are not   $\cD$-modules on $X$, potent categorical $G$-representation are not categories equipped with some form of $G$-action. 
Restriction
along the $S^1$-equivariant map $\wh \cL BG = \wh G/G \subset G/G = \cL BG$ gives a functor to  an {\em underlying} de Rham categorical representation
 $$
\Upsilon: G-cat^{pot} = 2\enhIndCoh(G/G)^{\Tate} \to 
2\enhIndCoh(\wh G/G)^{\Tate} \simeq  G-cat
 $$
but a potent categorical $G$-representation contains far more information.

Potent categorical $G$-representations $\cM$ also have 
a {\em global sections category} $\Gamma(\cM)$  or {\em category of  invariants} which naturally lifts to an object of 
$ 2\enhIndCoh(H\GIT W)$. 

\begin{ex}\label{constructing potent reps}
The main source of potent categorical $G$-representations 
 is given by linearizing $G$-stacks $X$. The rotation-equivariant map of loop spaces $\cL(X/G) \to \cL BG = G/G$ gives an object we denote by
${\mathbf D}^{pot}_G(X) \in G-cat^{pot}.
$
 Its  global sections and underlying de Rham categorical representation  are given respectively by potent $\cD$-modules
 and  $\cD$-modules 
  $$
   \Gamma({\mathbf D}^{pot}_G(X)) = \cD^{pot}(G\bs X)
   \qquad
   \Upsilon({\mathbf D}^{pot}_G(X)) =  \cD(X)$$

\end{ex}

\subsubsection{The zoo of circles}\label{zoo}
There are as many  variants on the notion of derived loop space $\cL X= \Map(S^1, X)$  as there are variants of the circle in the world of stacks. Here are some key examples of variants $A$ listed with their Picard (or 1-Cartier) duals $\Pic(A) =  \Hom(A, B\GG_m)$.
We also list our preferred  notation and name for the  mapping stack  $X^A = \Map(A,X)$.
 
\begin{itemize}
\item $A=S^1 = B\ZZ$, $\Pic (A) = \GG_m$, $X^{S^1} = LX$ loop space.
\item $A=B\GG_a$, $\Pic (A) = \wh\GG_a$, $X^{B\GG_a} =: L^u X$ unipotent loop space.
\item $A=B\wh\GG_a$, $\Pic (A) = \GG_a$, $X^{B\wh\GG_a} \simeq  T[-1]X$
graded loop space (odd tangent bundle).
\item $A=\SS^1_{dR} := (\GG_m)_{dR} = \GG_m/\wh\GG_a$, $X^{\SS^1_{dR}} =:L_{dR} X$ de Rham loop space.
\end{itemize}

\begin{ex}
For $X$ a scheme, all  the variants of loop spaces considered are canonically equivalent $ LX \simeq L_{dR} X \simeq L^uX \simeq TX[-1]$, and  they all coincide with the formal loop space  $\wh LX$. But for a stack they differ, though all naturally contain the formal loop space. For example, for $BG$, we find  
the moduli of connections $L_{dR} BG \simeq \Conn_G(\GG_m)$  and 
various adjoint quotients $LBG \simeq G/G$, $L^uBG \simeq G^\wedge_\cU/G$, where $\cU \subset G$ is the unipotent cone,  $T[-1] BG \simeq \frg/G$,  $\wh LBG \simeq \wh G/G$.
In particular, for $G=T$ a torus, we find
$L_{dR} BT \simeq \Conn_T(\GG_m)$, 
$LBT \simeq T \times BT$, $L^uBT \simeq T^\wedge_e \times BT$, where $e \in T$ is the identity,  $T[-1] BT \simeq \frt \times BT$,  $\wh LBT \simeq \wh T \times BT$.


\end{ex}

\begin{remark} For $BG$, we will focus on the 
 rotation-invariant open of regular singular connections $\Conn^{rs}_G(\GG_m) \subset \Conn_G(\GG_m)$
 within the de Rham loop space. Note the monodromy of connections gives an analytic equivalence from 
 $\Conn^{rs}_G(\GG_m)$ to the Betti loop space $G/G$.
 In particular, for $G=T$ a torus, we have 
 $\Conn^{rs}_T(\GG_m) \simeq \frt/\Lambda \times BT$.

\end{remark}

 For any variant of the circle $S^1 = B\ZZ$, one can implement a parallel notion of potent categorical $G$-representation. 
With this in mind, we sometimes refer to the default notion for  the Betti loop space $\Map(S^1,BG)=G/G$
 as {\em Betti} potent categorical $G$-representations
$$
G-{cat}^{pot} = G-{cat}^{pot,B}:= 2\enhIndCoh(G/G)^{\Tate}
$$
Here are two notable variants: 
For the de Rham stack $(\GG_m)_{dR}$, the ``de Rham loop space'' is the stack of connections 
$$\Map((\GG_m)_{dR}, BG) = \Conn_G(\GG_m) = \{t\partial_t + \frg[t,t^{-1}]\}/G[t, t^{-1}]$$ 
with rotation-invariant open of regular singular connections $\Conn^{rs}_G(\GG_m) \subset \Conn_G(\GG_m)$. We then have the notion of {\em de Rham}
potent categorical $G$-representations
$$
G-{cat}^{pot, dR} := 2\enhIndCoh(\Conn^{rs}_G(\GG_m) )^{\Tate}
$$
For the stack $B\wh \GG_a$, the ``graded loop space'' is the odd tangent bundle $\Map(B\wh \GG_a, X) = T[-1] X$, and we have the notion of {\em graded}
potent categorical $G$-representations
$$
G-{cat}^{pot, gr} := 2\enhIndCoh(\frg/G)^{\Tate}
$$
Note all of the above give extensions of the notion of de Rham categorical representations.


\subsection{Elaborations on $2\IndCoh$}\label{IndCoh elaborations}

We will need two key elaborations on $2\IndCoh$ -- {\em periodization} and {\em equviariant expansion} -- which we highlight here. 

\subsubsection{Periodization and Fourier transform}
As mentioned above, 
objects $\cM \in 2\IndCoh(X) $ have a singular support inside $T^*[2] X$. Thus one can microlocalize $2\IndCoh(X)$ to obtain a sheaf of 2-categories over $T^*[2] X$.
But we would like a richer version of $2\IndCoh(X)$ that lives 
over $T^* X$ and thus admitting objects corresponding to nonconic Lagrangians. 

To this end, we introduce the {\em periodic base} symmetric monoidal 2-category
$$
\cA  = 2\IndCoh^\diamondsuit (pt) 
:= 2\IndCoh(\GG_a)/2\QCoh(\GG_a)
$$
with convolution product given by addition in $\GG_a$. 
Then we introduce 
the {\em periodization} 
$$
2\IndCoh^\pi (X)=  2\IndCoh (X) \otimes  \cA
$$
regarded  as an $\cA$-module.

\begin{remark}
The periodization 
is an ind-coherent refinement, specifically in the context of  $2\IndCoh$, of the familiar periodization of
quasicoherent constructions by tensoring with $k[u,u^{-1}]$.
It is beyond our current scope to discuss further,
but strictly speaking, the  constructions involved 
should be performed inside $3\IndCoh(pt)$. 

Note the periodic base 
$\cA  = 2\IndCoh^\diamondsuit (pt) = 2\IndCoh(\GG_a)/2\QCoh(\GG_a)
$ is  the microlocalization of $2\IndCoh(\GG_a)$ away from the zero-section. Likewise,
 the periodization 
$
2\IndCoh^\pi (X)=  2\IndCoh (X) \otimes  \cA
$
is the parallel microlocalization of $2\IndCoh (X \times \GG_a)$. One can view it as living on the contact one-jet bundle $J^1X$. 
\end{remark}


A main theorem in the upcoming work~\cite{BZNS} is the following 2-Fourier transform which takes advantage of the above periodization.
As discussed below, a key application of this theorem, or in fact variations and extensions of it, will be in the analysis of potent categorical representations of dual tori.

\begin{theorem}\label{2-Fourier theorem}
For $V$ a finite-dimensional vector space  with dual $V^*$, there is a 2-Fourier transform
$$2\IndCoh^\pi(V)\simeq 2\IndCoh^\pi(V^*)$$
whose cyclic trace recovers a two-periodic version of the Fourier transform for $\cD$-modules
$\cD(V)\simeq \cD(V^*)$
\end{theorem}

\subsubsection{Equivariant expansion}\label{expansion section}
Recall for a scheme $X$,  one can view   $2\IndCoh(X)$ as an enlargement of $2\QCoh(X)$ with respect to microlocal parameters in the sense that $2\QCoh(X) \subset 2\IndCoh(X)$ consists of objects with singular support in the zero-section $X\subset T^*[2]X$.
As discussed in Remark~\ref{rem:2indcoh has many generators}, this reflects that even when $X$ is an affine scheme, 
objects of $2\IndCoh(X)$ are much more than a single category, but encode compatible collections of categories indexed by test schemes $Y\to X$. 

Similarly, for a stack $X$, we can introduce an enlargement
of $2\IndCoh(X)$ itself to a new 2-category $2\IndCoh^\epsilon(X)$ such that the inclusion 
$2\IndCoh(X) \subset 2\IndCoh^\epsilon(X)$  
is given by the completing  at the trivial equivariant parameter. The construction can be implemented by specifying that
objects of $2\IndCoh^\epsilon(X)$  encode compatible collections of categories indexed not only by test schemes but by select test stacks $Y\to X$  as well. 
In pursuing this,  we are guided by lessons from  equivariant homotopy theory as found in Elmendorf's theorem~\cite{Elmendorf}: to properly capture the geometry of group actions on spaces, one must pass from the more concrete (Borel) notion of equivariance as homotopy types with extra structure to considering diagrams of fixed-point spaces. 
The specific approach of capturing equivariance via pre-sheaves on test orbifolds appears explicitly in 
 global homotopy theory~\cite{henriquesgepner, jacobellipticIII, Schwedeglobal}.

For our current applications, we will in fact only need the concrete case of a quotient stack $X = Z/S^1$ where $Z$ is a scheme
with $S^1$-action (or similarly, the quotient by any of the zoo of circles collected in ~\S\ref{zoo}). In this case, to construct 
$2\IndCoh^\epsilon(Z/S^1)$, we 
 expand the generators of $2\IndCoh(Z/S^1)$ from schemes  $Y \to Z$ with $S^1$-action to include reductive abelian gerbes $BA \to Z$ with $S^1$-action. As a measure of the impact of this maneuver, in general $2\IndCoh^\epsilon(Z/S^1)$ is spread out  over the equivariant parameters $H^*(BS^1) \simeq k[u]$, while 
 $2\IndCoh(Z/S^1)$ is in fact $u$-torsion.

The following realization of 2-Cartier duality shows what equivariant expansion achieves in the base case  $X=BS^1$. It also shows how duality swaps 
equivariant expansion for microlocal expansion.
(On objects, it 
matches e.g. $B\GG_m/S^1 \to BS^1$, with $S^1$-action by $\mu\in \GG_m$, and $\mu^\BZ\to \GG_m$.)

\begin{prop}\label{expanded Mellin}
There is a canonical equivalence
$
2\IndCoh^\epsilon(BS^1) \simeq   2\IndCoh(\GG_m) 
$
restricting to an equivalence
$
2\IndCoh(BS^1) \simeq   2\QCoh(\GG_m),
$
and Tate localization corresponds to microlocalization
$$
2\IndCoh^\epsilon(BS^1)[u^{-1}] \simeq   2\IndCoh(\GG_m)/ 2\QCoh(\GG_m) 
$$
\end{prop}

\subsubsection{Summary}

Going forward, we will write $2\enhIndCoh$ for the result of both the periodization $2\IndCoh^\pi$ and the equivariant expansion
$2\IndCoh^\epsilon$. In our developments, the latter is only relevant when we speak about loop rotation equivariance and can  be omitted otherwise. 


\subsection{The potent Hecke monad}\label{Hecke monad section}

Given the above elaborations $2\enhIndCoh$  on $2\IndCoh$ -- encompassing the periodization  $2\IndCoh^\pi$ and the  equivariant expansion $2\IndCoh^\epsilon $ -- we can now  begin to analyze Betti, de Rham, and graded potent categorical representations 

$$
G-{cat}^{pot} = G-{cat}^{pot, B} := 2\enhIndCoh(G/G)^{\Tate}
$$
$$
G-{cat}^{pot, dR} := 2\enhIndCoh(\Conn^{rs}_G(\GG_m) )^{\Tate}
\qquad
G-{cat}^{pot, gr} := 2\enhIndCoh(\frg/G)^{\Tate}
$$
using paradigms of representation theory, specifically induced representations and  Hecke algebras.

\subsubsection{Potent highest weight theory}

Traditional representation theory of a reductive group $G$ is stratified by parabolic induction from Levis, starting with principal series induced from tori and culminating in cuspidal representations.  
Geometric representation theory 
on the other hand ``lives'' on the flag variety  (cf.~Deligne-Lusztig theory, character sheaves, Beilinson-Bernstein localization)
which is ``home" to intertwiners on principal series.
This is reflected in categorical representation theory by the fact that 
all de Rham $G$-categories appear in principal series~\cite{BZGO}: there is a Morita equivalence between the de Rham group algebra $\cD(G)$ and the universal  Hecke category $\cH^{univ}_G = \cD(U\bs G/ U)$ of intertwiners on principal series $\cD(G/U)$.

This phenomena extends to potent categorical representations of all types, but for simplicity, we will focus here on the default Betti type. 
Passing to  loops in  the fundamental correspondence $BG \leftarrow BB \to BT$, 
provides the following:

\begin{defn}
The {\em Betti potent  Hecke monad} 
is the  algebra object  $$\cH_G = \cH_G^B \in 2\enhIndCoh(T/T \times T/T)^\Tate \simeq  T\cat^{pot}\otimes T\cat^{pot}$$
represented by the rotation-equivariant groupoid 
$$
\cL(B\bs G/B) \to  \cL(BT \times BT) = T/T \times T/T
$$
or equivalently, in the notation of \S\ref{constructing potent reps},
 by  the potent representation
$${\mathbf D}^{pot}_{T\times T}(U\bs G/U) \in T\cat^{pot}\otimes T\cat^{pot}$$

The {\em de Rham} and {\em graded} potent  Hecke monads
$\cH_G^{dR}$ and $\cH_G^{gr}$ are defined accordingly.

\end{defn}

The potent  Hecke monad $\cH_G$ is not itself a category, but knows about many  Hecke categories. More specifically, completing along a Lagrangian (taking morphisms to $\cH_G$ from  objects of  
$2\enhIndCoh(T/T \times T/T)^{\Tate}$) recovers many  
familiar Hecke categories. 
Notably, 
taking morphisms from the object associated to
${e}/T \times {e}/T$ gives the universal  Hecke category $\cH^{univ}_G = \cD(U\bs G/U)$, a $\cD$-module version 
 of the universal monodromic Hecke category of~\cite{tayloruniversal}. 
Note the bimodule structure over $\End(\cM({e}/T)) \simeq\cD(T) \simeq \QCoh(\check\frt /\check\Lambda)$ allows one then to further specialize the monodromicity. In particular, one can recover
the pro-unipotent Hecke category $\cH^{unip}_G$.
For another example, taking morphisms from the object associated to $T/T \times T/T$, 
i.e.~global sections,
gives the potent Hecke category $\cH^{pot}_G = \cD^{pot}(B\bs G/B)$, 
a cyclic homology version of the $K$-motive Hecke category of~\cite{eberhardtK,eberhardteteve}. Completing along equivariant parameters gives the equivariant Hecke category $\cH^{eq}_G$.


Now the Morita theorem for categorical de Rham representations~\cite{BZGO} extends to the following potent statement:

\begin{theorem}\label{thm highest weight} Via loops in the fundamental correspondence $G/G \leftarrow B/B \to T/T$, potent categorical $G$-representations
$G-cat^{pot} = 2\enhIndCoh(G/G)^{\Tate}$
are monadic over potent categorical $T$-representations $T-cat^{pot}  = 2\enhIndCoh(T/T)^{\Tate}
$.
The monad is the potent  Hecke monad $\cH_G$.

A similar statement holds for de Rham and graded 
potent categorical representations.
\end{theorem}

\subsection{Potent Langlands duality}\label{potent conjecture section}
Here we discuss results and conjectures on a Langlands correspondence for potent categorical representations.
We start with the case of tori.

\subsubsection{Multiplicative Fourier transform}

Recall the 2-Fourier transform 
of Theorem~\ref{2-Fourier theorem}  categorifies a version of the Fourier transform for $\cD$-modules.
The following multiplicative version is analogous instead to the $q$-difference Fourier transform.

For a torus $T$  with Lie algebra $\frt$ and coweight lattice $\Lambda$, note that regular singular $T$-connections take the concrete form $\Conn_T^{rs}(\GG_m) \simeq \frt/\Lambda \times  BT$.

\begin{theorem}\label{difference Fourier}
For a torus $T$, there is a multiplicative 2-Fourier transform
$$T\cat^{pot, dR}=
2\enhIndCoh (\frt/\Lambda \times  BT)^{\Tate} \simeq 
2\enhIndCoh(\check\frt/\check \Lambda \times B\check T)^{\Tate} 
=
\Tv\cat^{pot, dR}$$ 
\end{theorem}

\begin{remark}
One can view the theorem as the concatenation of two kinds of results: first,  duality equivalences such as
$$
2\enhIndCoh(\frt) \simeq 
2\enhIndCoh(\check \frt)
\qquad
2\enhIndCoh (BT) \simeq 
2\enhIndCoh(B \Lambda_{\check T}) 
$$
and then calculations of the Tate constructions appearing. 
We expect with sufficient analytic foundations  there will  be an  analogous equivalence  in the  Betti setting.
\end{remark}

\begin{remark}\label{rem:gr t duality}
There is  also an analogous equivalence for graded potent categorical $T$-representations
$$
T-cat^{pot, gr} = 2\enhIndCoh(\frt/T )^{\Tate} \simeq 
2\enhIndCoh(\check\frt/\check \Lambda)
$$
whose cyclic trace recovers a version of  
the Mellin transform. See Proposition~\ref{expanded Mellin} for a related statement: up to algebraic differences, specifically 
between $\frtv/\Lambda$ and $\Tv$, graded potent categorical $T$-representations are equivalent to  categorical $T_c$-representations as formalized by $2\enhIndCoh(BT_c)$.

\end{remark}

\begin{ex}
Within the theorem, we find the familiar equivalence for de Rham categorical $T$-representations 
$$
T-cat = 2\enhIndCoh(\wh T /T)^{\Tate} \simeq  2\QCoh(\check \frt/\check\Lambda)
$$
deduced from the monoidal equivalence $\cD(T) \simeq \QCoh(\check \frt/\check\Lambda)$.
On the right hand side,  the parameters of the theory appear clearly as the positions $\check \frt/\check\Lambda$. These match with the characters of de Rham categorical representations of $T$
on the left hand side.
But the situation is  asymmetric:  the parameters of the theory do not include the equivariant parameters/momenta $\frt \simeq \check \frt^* $.
\end{ex}

\subsubsection{Main conjecture}
For a general reductive group, we can invoke the highest weight theory of Theorem~\ref{thm highest weight} and the multiplicative 2-Fourier transform of Theorem~\ref{difference Fourier} to arrive at the following.

\begin{conj}\label{potent duality conj} 
For a reductive group $G$, the  multiplicative 2-Fourier transform of Theorem~\ref{difference Fourier}  identifies the de Rham potent Hecke monads $\cH_G^{dR}$ and $\cH_{\check G}^{dR}$. Consequently, by Theorem~\ref{thm highest weight}, there is an equivalence of de Rham potent categorical representations  
$
G-cat^{pot, dR}\simeq \Gv-cat^{pot, dR}.
$

\end{conj}

\begin{remark}
We expect with sufficient analytic foundations  there will be an  analogous equivalence  in the  Betti setting. In the graded setting, we expect a relation to topological categorical representations as discussed below in \S\ref{Coulomb section}.   
\end{remark}

The conjecture says that there exists a {\em single} Hecke monad living over the symplectic space $\frt/\Lambda \times \check\frt/\check \Lambda $ from which  different versions of  Hecke categories (universal, $K$-theoretic, monodromic, equivariant, etc.) arise by pairing with some Lagrangian and then prescribing parameters. Moreover, each of these different versions will have two Langlands dual realizations depending on how we parse them. 

For example, on the one hand,
pairing $\cH_G^{dR}$ with
$0 \times \check \frt/\check \Lambda$ gives the universal  Hecke category $\cH^{univ}_G = \cD(U\bs G/U)$, a $\cD$-module version 
 of the universal monodromic Hecke category of~\cite{tayloruniversal}. 
 On the other hand, 
 pairing $\cH_{\check G}^{dR}$ with
$0 \times \check \frt/\check \Lambda$ gives the
 potent Hecke category $\cH^{pot}_{\check G} = \cD^{pot}(\check B\bs \check G/\check B)$, 
a Hochschild homology version of the $K$-motive Hecke category of~\cite{eberhardtK,eberhardteteve}.
The conjecture thus predicts an equivalence between these two, a variation on the universal Koszul equivalence of~\cite{tayloruniversal}. 

\section{Inspiration from topological field theory}\label{TFT section}


In this section, we discuss  some of the physics motivating the developments of the prior section. 
Our working context is  fully extended 3d topological quantum field theories (TQFTs) and their dualities, specifically the nonabelian duality of  {\em 3d mirror symmetry} intermediate between 2d mirror symmetry and 4d  electric-magnetic/Langlands duality. 

\subsection{Ramping up: from 2d to 3d}
To orient the discussion, we  recall in briefest terms   2d  mirror symmetry and the structure of fully extended TQFTs via the Cobordism Hypothesis.


\subsubsection{2d mirror symmetry}
Homological  (or fully extended 2d) 
mirror symmetry, as envisioned by Kontsevich~\cite{KontsevichICM}, takes the form of equivalences between categories of A-branes (Fukaya categories, $\cD$-modules, microlocal sheaves,...) on symplectic manifolds and B-branes (coherent sheaves, matrix factorizations,...) on mirror complex manifolds. 
Such branes form the local boundary conditions or {\em boundary theories}
in 2d TQFTs, and
the corresponding categories organize their
 interfaces.
 Following~\cite{costello, KontsevichSoibelman}, 
 the entire structure of a 2d TQFT is determined by its  category of branes,
 and
 simple conditions -- smoothness and properness (for framed theories) and Calabi-Yau structure (for oriented theories) -- are sufficient to qualify a category to be the category of 
branes in a  2d TQFT.  

From this general perspective, one can view 2d mirror symmetry as the challenge: given a 2d TQFT of any source, find a B-model equivalent to it.
Broadly speaking, there are two natural sources of 2d TQFTs one would like to  tackle in this way: 
$\sigma$-models 
and gauge theories. 
For $\sigma$-models, one can often express the boundary theories as sheaves of some kind on a target;  for gauge theories,  the
boundary theories take the form of representations of a group.
Of course, the language of  
gauged $\sigma$-models and technology of stacks provides a uniform approach to the two, so that ``everything is a $\sigma$-model" but with a possibly non-traditional target.

For our developments, the following is a small but key example
of 2d mirror symmetry relating a gauge theory to a $B$-model.

\begin{ex}  \label{ex basic 1}
     The branes in 2d A-twisted supersymmetric $U(1)$-gauge theory  are  ``topological $U(1)$-representations": (bounded) cochain complexes with an action of $U(1)$, i.e.~an action of chains $C_{-*}(U(1)) \simeq k[v], |v|=-1$, as arise for example by taking cochains $C^*(X)$ on $U(1)$-spaces $X$.

     The mirror B-model has target the graded affine line $\AA^1[2]$ with functions given by  equivariant cohomology $H^*(BU(1)) \simeq k[u]$, $|u|=2$. Its branes are coherent sheaves on $\AA^1[2]$,
     i.e.~(bounded)  complexes of $k[u]$-modules.
          
The mirror equivalence between the  categories of branes is given by Koszul duality: taking the $U(1)$-invariants of a brane on the A-side gives a brane on the B-side.

\end{ex}

\subsubsection{Fully extended TQFTs}
The powerful tools of sheaf theory -- such as functoriality and descent -- have played a vital role in the study of 2d TQFTs and
their mirror symmetry. 
We would like to study 3d TQFTs  with analogous tools. 

To start, thanks to the Cobordism Hypothesis~\cite{BaezDolan,jacobTFT}, 
a fully extended $d$-dimensional TQFT $Z$ is equivalent data to the knowledge of the higher categorical object $Z(pt)$, which organizes its collection of boundary theories and their interfaces. In other words, we can think of a $d$-dimensional TQFT informally as the study of a collection of $(d-1)$-dimensional TQFTs with extra structure comprising its boundary theories. 

From this perspective, the challenge in 3d is to construct physically meaningful 3d TQFTs by prescribing their 2-categories of boundary theories -- collections of 2d TQFTs with extra structure --  and then realize predicted dualities by finding equivalences among these 2-categories.
As with 2d TQFT, 
there are broadly speaking two natural sources of 3d TQFTs: 
$\sigma$-models, whose boundary theories are given by {\em geometric families of 2d TQFTs}, and
 gauge theories,
whose   boundary theories are  given by {\em 2d TQFTs with symmetry}, i.e.~categorical representations. In both cases, the  precise flavor of these constructions will depend on twists, with for example de Rham categorical representations arising in A-twisted gauge theory.



 Higher sheaf theory offers a promising extension of sheaf theory
  to study fully extended 
 theories in all dimensions, in particular to construct  
higher categorical objects organizing
 boundary theories and interfaces. 
In the case of 3d TQFTs, we will sketch how the notion of {\em sheaves of categories} achieves this in specific cases of interest.

\subsection{3d TQFT via 2-categories}\label{3d tqft via shv}

In this section, we reinterpret the constructions of ~\S\ref{sect:pcrs} and ~\S\ref{IndCoh elaborations}
to show how  {\em sheaves of categories} and {\em categorical representations}
 provide mathematical models of {\em Rozansky-Witten theory} and {\em 3d gauge theory}.

{\color{red} 

}


\subsubsection{Rozansky-Witten theory}\label{RW theory}

Given a holomorphic symplectic manifold $M$, there is a $\ZZ/2$-graded 3d TQFT $RW_M$, called  {\em Rozansky-Witten (RW) theory} with target $M$, which is a topologically twisted form of the 3d hyperk\"ahler $\sigma$-model. Much is known about RW theory including the key structures:
\begin{itemize}
\item functions on $M$ are the local operators: $RW_M(S^2)\simeq \cO(M)$

    \item the B-model of $M$ is its compactification on the circle: $RW_M(S^1)=\QCoh(M)$
\end{itemize}

 Kapustin, Rozansky and Saulina~\cite{KRS,KR, KapustinICM}) pioneered the study of RW theory as a {\em fully extended} TQFT by describing boundary theories, objects of the putative 2-category $RW(M):=RW_M(pt)$ assigned to a point. 
 Their proposed objects are (decorated) Lagrangians  $L\subset M$, and morphisms $L_1 \to L_2$, are categories of matrix factorizations (when $L_2$ is graphical over $L_1$). 
In particular, for a single Lagrangian $L \subset M$, one finds its endomorphisms are $\QCoh(L)$. Hence the completion of the theory along $L$ is given by $2\QCoh(L) = \QCoh(L)-\module$.

Thus for a cotangent bundle  $M= T^*L$, a reasonable first guess at $RW(T^*L)$ is simply quasicoherent sheaves of categories on the zero-section $2\QCoh(L)$ (see also~\cite{TVChern,BZFN}).
However, the compactification of $2\QCoh(L)$  on the circle is not $\QCoh(T^*L)$ but its completion $\QCoh_L(T^*L)$ consisting of sheaves supported along the zero-section. So this completion of RW theory is missing the momenta parameters and 
thereby missing all of the
branes supported away from the zero-section. One needs an  expansion of $2\QCoh(L)$ to correct for this completion.



Among many structures, Stefanich's higher sheaf theory~\cite{stefanichIndCoh}
assigns a 2-category $2\IndCoh(L)$  
as an expansion of $2\QCoh(L)$.
It enjoys the expected properties of $RW(T^*L)$, including the correct local operators and compactification on the circle (and indeed its 2-periodicity).
As mentioned in~\S\ref{sect:pcrs} and explained in~\S\ref{IndCoh elaborations},
a periodized variant $2\IndCoh^\diamondsuit(L)$ enjoys the additional property that Lagrangians in $T^*L$ define objects
(as well as functoriality for Lagrangian correspondences, for example associated to proper maps from smooth varieties to $L$).
and 
that morphisms between graphical Lagrangians can be calculated by matrix factorizations~\cite{Yiu}. Furthermore, the 2-Fourier transform of  Theorem~\ref{2-Fourier theorem} is a first instance of the expected symplectic invariance of 
 $2\IndCoh^\diamondsuit(L)$ as a model for $RW(T^*L)$.

\begin{ansatz}\label{rw ansatz}
The boundary theories for RW theory with target a cotangent bundle 
are periodized ind-coherent sheaves of categories on the base:
$
RW(T^*L) = 2\IndCoh^\diamondsuit(L).
$
\end{ansatz}

\begin{remark}\label{rem:gauged ansatz}
The Ansatz naturally extends to gauged RW theory: if $L$ is a $G$-variety,  then we model  $G$-gauged RW theory of $L$  by 
periodized ind-coherent sheaves of categories on the quotient stack:
$
RW^G(T^*L) = 2\IndCoh^\diamondsuit(L/G).
$
Of notable interest for us is the $G$-gauged RW theory of $T^*G$ for the adjoint $G$-action on $G$.
\end{remark}

Going further, elaborations on the above constructions  provide a model
for RW theory with target any open $\Omega \subset T^*L$ as follows.
The fact that $2\IndCoh^\diamondsuit$ has local operators $\cO(T^*L)$ allows us to microlocalize $2\IndCoh^\diamondsuit$ over 
$T^*L$.  Hence, inspired by Kashiwara-Schapira's microlocal theory of  sheaves on manifolds~\cite{KashiwaraSchapira}, we arrive  at the following proposal:
$$
RW(\Omega) = 2\IndCoh^\diamondsuit(L)/\{ \text{branes supported away from } \Omega \}  
$$
We conjecture this is a symplectic invariant suitably understood.

Going further still, one can consider symplectic manifolds $M$ that 
may  be obtained from an open $\Omega \subset T^*L$ by Hamiltonian reduction.   Here the functoriality of $ 2\IndCoh^\diamondsuit$ allows one to  arrive at a monadic definition for $RW(M)$ in terms of $RW(\Omega)$. We conjecture this too will be a symplectic invariant suitably understood. Finally, turning on equivariant parameters for topological circle actions using the equivariant expansion of $2\IndCoh$ allows us to construct more TQFTs as {\em mass deformations} of RW theories, see \S~\ref{Omega section}.

\subsubsection{3d gauge theory}

A central theme of representation theory -- as realized by geometric quantization and the orbit method -- is that representations of groups arise via symmetries of quantum mechanical systems, i.e.,~1d quantum field theories. 
In direct analogy, 2d quantum field theory provides a primary source of categorical representations: the starring role of Hilbert spaces of states in quantum mechanics is now played by categories of boundary theories; when a 2d theory carries global symmetry (for example, as a $\sigma$-model into a space with symmetry), its category
of boundary theories becomes a categorical representation of some flavor.


 
%


Here we enter the realm of 3d TQFT when we organize all  categorical representations of a given flavor into the 2-category of boundary theories of a 3d gauge theory. Teleman's pioneering work~\cite{TelemanICM,TelemanCoulomb} focused on A-models of symplectic manifolds with Hamiltonian symmetry for a compact group~$G_c$. Their Fukaya categories are topological $G_c$-categories and provide important examples of boundary theories for twisted 3d $\cN=4$ Yang-Mills theory. Here a topological $G_c$-category means a module category for the convolution monoidal category  $\Loc(G_c)$ of local systems on $G_c$. 
However, Teleman crucially discovered that the 2-category of $\Loc(G_c)$-modules is missing half the parameters, specifically the equivariant parameters, as we revisit from the 3d mirror picture in \S~\ref{Coulomb section}.

The theory of de Rham categorical representations, i.e. modules for $\cD(G)$, provides boundary theories for a different gauge theory: a topological twist of the 3d maximally supersymmetric ($\cN=8$) Yang-Mills theory (the $\cN=4$ gauge theory with adjoint matter). But again this direct approach misses half the parameters, whence the need for the potent categorical representations introduced in
~\S\ref{sect:pcrs} and ~\S\ref{IndCoh elaborations}. We propose the graded variant $G-cat^{pot,gr}$ as the full 2-category of boundary theories in this extended TQFT (while $G-cat^{pot}$ itself models the equivariant compactification of 4d $\cN=4$ Yang-Mills on a circle, see ~\S\ref{4d section}).
The key source of such boundary theories come from $G$-spaces $X$ via the construction
${\mathbf D}^{pot}_G(X) \in G-cat^{pot}.
$
Informally speaking, this is a family of 2d TQFTs over the equivariant parameters $T\gitquot W$ whose fibers are categories of $\cD$-modules on fixed-point stacks, and whose global sections is the category of potent $\cD$-modules.

 \subsection{3d mirror symmetry}\label{3d mirrors}
 Now we turn to 3d mirror symmetry~\cite{IntrilligatorSeiberg} to begin to give the physical context for
  the potent duality of Conjecture~\ref{potent duality conj}. Inspired by 2d mirror symmetry, we take the viewpoint: given a 3d TQFT $Z$, our goal for 3d mirror symmetry is to find a holomorphic symplectic target $M$ so that $Z$ is equivalent to the Rozansky-Witten theory $RW_M$ as constructed in ~\S\ref{RW theory}. This perspective on 3d mirror symmetry as a fully-fledged higher analogue of homological mirror symmetry,
i.e.~as equivalences of 3d fully extended TQFTs, emerged in~\cite{TelemanICM,BDGH, hilburngammagemazelgee,hilburngammagehypertoric,pascaleffsibilla}. It provides a rich physical context for the fundamental phenomenon of symplectic duality in geometric representation theory~\cite{BLPW}.



 \subsubsection{Moduli of vacua} 
Given a 3d TQFT $Z$, to find a mirror holomorphic symplectic target $M$, we begin with an important lesson from physics:
a 3d TQFT $Z$ comes  with an intrinsic parameter space $\cM_Z$, called its {\em moduli  of vacua},  that is an affine graded Poisson  variety. 
By definition, the coordinate ring of $\cM_Z$ is given by
(the cohomology of) the ring of local operators $\cO(\cM_\cZ)=\cZ(S^2)$. The graded Poisson structure on $\cM_Z$ comes from the $E_3$-structure on the local operators. 

 The moduli of vacua $\cM_Z$, viewed as the target of a $\sigma$-model,  provides a ``low-energy approximation/affinization" of the original theory $Z$, in the sense
 that we expect $Z$ to look like 
 RW theory with target $\cM_Z$ or a resolution of it.

 \subsubsection{Coulomb branches} \label{Coulomb section}
For A-twisted supersymmetric $\cN=4$ 3d gauge theories $Z$, the moduli of vacua $\cM_{Z}$ is the {\em Coulomb branch}, for which an extremely influential mathematical construction was given by Braverman, Finkelberg and Nakajima~\cite{BFN}. For $G$-gauge theories, the Coulomb branch is birational  to $T^*\check T\gitquot W = (\Tv \times \frt)\GIT W$. For equivariant compactifications of 4d $G$-gauge theories, as discussed in~\S\ref{4d section}, one  encounters 
 the ``multiplicative/K-theoretic'' Coulomb branch which is birational  to 
 $(\check T \times T)\gitquot W$. 
 
In the case of pure gauge theory, the Coulomb branch is identified (thanks to~\cite{BFM}) with Kostant's completed phase space $J_{\Gv}$ of the Toda lattice for the dual group~$\Gv$. Teleman~\cite{TelemanICM} suggested that the full theory of topological $G$-categories should thus be equivalent to the Rozansky-Witten theory of $J_{\Gv}$. However he discovered that $\Loc(G_c)$-modules only see the completion of $RW(J_{\Gv})$ along a Lagrangian, whence the need for expansion.



 \begin{ex}\label{ex: 3d u(1)} 
We illustrate how the genuine version of $2\IndCoh$ introduced in~\S\ref{expansion section} realizes the expansion of $A$-twisted $U(1)$-gauge theory envisioned in~\cite[\S 4]{TelemanICM}.
 We take as the boundary theories $Z(pt)$ categories with a topological $U(1)$-action
 in the genuine sense modeled by 
 $
2\enhIndCoh(BU(1)).
$
Then one finds (as in~\cite{TelemanICM})
$$
Z(S^2) = C_{-*}(S^1\bs L S^1/S^1) \simeq k[z, z^{-1}, u],
\quad |z| = 0, |u| = 2
$$
so that $\cM_Z = T^*\GG_m = J_{\GG_m}$.
Following
 Proposition~\ref{expanded Mellin},
 we indeed have 
$$
Z(pt) = 2\enhIndCoh(BU(1)) \simeq   2\enhIndCoh(\GG_m) = RW(T^*\GG_m)
$$ 

Note   $\Loc(U(1)) \simeq \QCoh(\GG_m)$, and hence $\Loc(U(1))$-modules only see the completion 
$2\QCoh(\GG_M) \subset 2\enhIndCoh(\GG_m) = RW(T^*\GG_m)$ along the zero-section $\GG_m \subset T^*\GG_m = J_{\GG_m}$.
 \end{ex}

%

We now informally (and analytically) discuss a proposed mirror symmetry for the graded potent theories $G-cat^{pot,gr}$, modeling 3d $\cN=8$ gauge theory. (Its multiplicative version, the proposed duality for $G-cat^{pot}$ of Conjecture~\ref{potent duality conj}, is discussed below in \S\ref{4d section}.)   

For an abelian gauge group $T$,
as discussed in Remark~\ref{rem:gr t duality},
the  mirror equivalence with RW theory of the Coulomb branch holds:
$
T-cat^{pot, gr} \simeq 
 RW(T^* \check T)
$ (and in fact agrees with Teleman's expectations for mirror symmetry for topological $T_c$-categories as in Example~\ref{ex: 3d u(1)}).
For general gauge group $G$, the Coulomb branch for the $\cN=8$ theory is 
$T^*\check T\gitquot W = (\Tv \times \frt)\GIT W$. Thus we should look for mirror symmetry to relate $G-cat^{pot,gr}$ with a suitable model $RW^{\check G} (T^* \check G)^{0}$ of the RW theory for this space. Here the superscript $0$ reflects that we take the ``zero-mode" of the $\check G$-gauged RW theory, 
a limit at zero of the mass deformation discussed below in \S\ref{Omega section} in analogy with  the notion of ``limit category" coming from DT theory ~\cite{PadurariuToda}. 

Within this conjectural equivalence, 
we also expect one can extract Teleman's proposed equivalence for 
expanded topological $G_c$-categories: on the $A$-side, the affinization projection $\frg /G\to \frt \gitquot W$ (passing to global sections as in \S\ref{def of potent reps section}) leads to a realization of topological $G_c$-categories as a monadic shadow of $G-cat^{pot,gr}$; on the $B$-side, this corresponds to the monadic shadow of $RW^{\check G} (T^* \check G)^{0}$ given by restriction to the regular locus 
$T^*_{\reg}\check G \subset T^*\check G$. In fact, this matching reflects an expanded ``Whittaker normalization" of the duality for potent Hecke monads found in Conjecture~\ref{potent duality conj}.

 \subsection{Tate, mass and $\Omega$}\label{Omega section}
A key role throughout the developments of \S\ref{GRT intro} is played by the Tate construction, i.e., the process of turning on equivariant parameters for circle actions. This process appears naturally in TQFT for both external (global) and internal (gravitational) symmetries, as {\em mass deformations} and the $\Omega$-{\em background}.

First, suppose we have a TQFT with a topological action of $S^1$, which we can encode by an $S^1$-action on its collection of boundary theories $\cC$. The $S^1$-invariants $\cC^{S^1}$ form a family over the $u$-line where $u$ is the equivariant parameter $H^*(BS^1) \simeq k[u], |u|= 2$. This family encodes a family of TQFTs, a {\em mass deformation} of the original theory. 
The Tate construction results from giving $u$ a non-zero value, realized algebraically by inverting $u$. This often has the effect of localizing the theory to fixed points (i.e., making some of the fields massive). 

The Nekrasov $\Omega$-background~\cite{Nekrasov,NekrasovWitten} provides a mechanism to turn on equivariant parameters for symmetries of the spacetime itself. This has the remarkable feature of deforming $B$-type twists into $A$-type twists. For 3d TQFTs, the $\Omega$-background on $S^1$ (or $S^2$) results in the deformation quantization of the category of sheaves (or ring of functions) on the moduli space of vacua~\cite{Yagi}. For RW theory on $M=T^*X$ this precisely recovers the Koszul duality between $S^1$-equivariant coherent sheaves on loop spaces ($B$-type) and $\cD$-modules ($A$-type) of \S~\ref{loops and conns section}.   (For stacks $X$ -- i.e., for {\em gauged} RW theory -- we find instead a deformation of sheaves on the inertia stack $\cL X$ to the potent $\cD$-modules of \S~\ref{loops and conns section}.)
The resulting quantized Coulomb branch construction recovers many algebras of interest, including Cherednik algebras in the case of gauge theory with adjoint matter~\cite{BLPW,BFN}.

We can also view this ``internal'' $S^1$-action for a 3d TQFT on $S^1$ as ``external'', i.e.~as a global symmetry of the compactification of $RW_M$ on $S^1$, the 2d B-model of $M$. Turning on the equivariant parameter then deforms this B-model to an A-model. One dimension higher, our source of A-type 3d gauge theories is to start with a 4d B-model, compactify on $S^1$ and turn on mass deformation for the rotational symmetry. We built our model $G-cat^{pot, gr}$ for topologically twisted $\cN=8$ gauge theory as a mass deformation (Tate construction) of the $G$-gauged RW theory into $T^*(\frg/G)$, turning the 3d B-model into an $A$-twisted gauge theory. Here the $S^1$-symmetry came from realizing $\frg/G$ as a (graded) loop space into $pt/G$, i.e.~realizing the RW theory $RW_{T^*(\frg/G)}$ as a graded form of the compactification on $S^1$ of a 4d gauge theory $\cB_G$ (see \S~\ref{4d section} below). Our construction of the theory of potent $G$-categories $G-cat^{pot}$ is likewise  a mass deformation of the gauged RW theory for the cotangent to the (Betti) loop space $\cL(BG)\simeq G/G$, a 3d shadow of a 4d B-model in $\Omega$-background. If we then consider this 3d theory itself in an $\Omega$-background on $S^1$ (thus the 4d theory on an equivariant $S^1\times S^1$), we expect to encounter the category of modules for the quantized K-theoretic Coulomb algebra for gauge theory with adjoint matter, which is (a specialization of) the double affine Hecke algebra~\cite{BFM,branesDAHA}, see \S~\ref{DAHA appears}.

\subsection{Equivariant compactification in 4d}\label{4d section}

 Finally, 
we discuss here the physical origin of the potent Langlands duality of Conjecture~\ref{potent duality conj}.



\subsubsection{Geometric Langlands, very briefly}
Much of our thinking on  geometric representation theory has been influenced by Kapustin-Witten's discovery~\cite{KapustinWitten} of a physical origin for the geometric Langlands correspondence, specifically as an aspect of electric-magnetic or S-duality for 4d supersymmetric gauge theory. The maximally supersymmetric ($\cN=4$) 4d Yang-Mills theory has a topological shadow that is a 4d TQFT depending on a parameter $\Psi\in \PP^1$. Within this family are the $A$-twist $\cA_G$ ($\Psi=0$) and   $B$-twist $\cB_G$ ($\Psi=\infty$),  which are 4d analogues of the two sides of mirror symmetry in 2d and symplectic duality in 3d. The S-duality conjecture (at $\Psi=0$) is an equivalence of 4d TQFTs $\cA_G\simeq \cB_{\Gv}$ for Langlands dual groups. 
We proposed
the Betti variant of geometric Langlands~\cite{BettiLanglands}
to try to mathematically model some of the physics Kapustin-Witten describe.

For a surface $\Sigma$, the equivalence $\cA_G(\Sigma)\simeq \cB_{\Gv}(\Sigma)$ (along with many variants of it) is now a landmark theorem of~\cite{GLCproof}.
For the circle $S^1$, 
the expected equivalence $\cA_G(S^1)\simeq \cB_{\Gv}(S^1)$ is 
the {\em local Geometric Langlands conjecture} (expanded so as to drop the nilpotent singular support condition, cf.~\cite[Remark 3.16]{BettiLanglands}). The A-side
$\cA_G(S^1)$ is a 2-category of categorical $LG$-representations -- which we expect can be modeled by a potent version of module categories for the affine Hecke category. The B-side $\cB_\Gv(S^1)$ is the gauged RW theory $RW^{\check G}(\check G)$ -- which following Ansatz~\ref{rw ansatz} can be modeled by the 2-category
$2\enhIndCoh(\Gv/\Gv)$ 
of ind-coherent sheaves of categories on $\Gv$-local systems on $S^1$.
(In fact, Stefanich developed the theory of higher ind-coherent sheaves of categories to describe the entire $B$-model $\cB_\Gv$ as the fully extended TQFT with value on a point $\cB_\Gv(pt)=3\IndCoh(B\Gv).$)

\subsubsection{Back to 3d, equivariantly} A theme of our research over the years
has been to start with the expected equivalence $\cA_G(S^1)\simeq \cB_{\Gv}(S^1)$ and to explore the outcome
of turning on rotation equivariance.
As a starting point, we observed~\cite{conns}
that Bezrukavnikov's theorem ~\cite{romahecke} on affine Hecke categories
-- governing the tamely ramified part of the local geometric Langlands correspondence -- 
becomes the Koszul duality of finite Hecke categories~\cite{BGS,BezYun}. 
Key to this
 was the mechanism summarized in~\S\ref{loops and conns section}:
turning on the $\Omega$-background deformed B-side categories of branes into categories of  A-branes, leading to  symmetric equivalences between A-models on both sides.
The developments of this paper arose from trying to understand Witten's discovery~\cite{WittenAtiyah} of a physical explanation for this phenomenon. 
Witten showed using supersymmetry that after $S^1$-equivariant compactification, i.e. compactifying on $S^1$ and turning on the $\Omega$-background, the Kapustin-Witten theories
for $\Psi\in \PP^1$ all become identified, including the
 $A$-twist $\cA_G$ ($\Psi=0$) and   $B$-twist $\cB_G$ ($\Psi = \infty$):
 $\cA_G(S^1)^{\Tate}\simeq \cB_G(S^1)^{\Tate}.
 $ 
 When combined with S-duality, this results in a {\em symmetric duality}, which one can interpret as identifying A-twisted 3d theories for Langlands dual groups. Unwinding this, we arrived at the 3d
duality of Conjecture~\ref{potent duality conj} for potent categorical representations. 
From our perspective, a remarkable aspect of this is the
key guidance 4d provides about 3d, both in predicting a Langlands duality in the first place and in the natural appearance of 
elusive equivariant parameters (as captured by equivariant cohomology or $K$-theory) under equivariant compactification.

\section{Potent representation theory}\label{applications}

In this section, we discuss some possible impacts of potent categorical representations in representation theory, organizing the discussion with parallels with the theory of Hecke algebras.
We formulate conjectures at three different levels of categoricity corresponding to viewing potent categorical representations as a 3d  TQFT: module categories (boundary conditions associated to the point), trace/cocenter categories (associated to the circle), and homology of moduli spaces (associated  to surfaces). 

Much of the discussion of this section is informal, and following the conventions stated in \S\ref{conventions},
we will implicitly work two-periodically and not always spell out precise algebraic forms of conjectures, settling for analytic forms, for example eliding the  distinction between $T$ and $\frt/\Lambda$. 



\subsection{Finite, affine, and double affine Hecke algebras}
In this preliminary section, 
we orient the discussion by
locating potent categorical representations 
within the three flavors of Hecke algebras associated to a reductive group. 
To do so, we introduce the following table then  make some clarifying comments about its contents.

\renewcommand{\arraystretch}{1.5}
\begin{table}[h!]
    \centering

    \begin{tabular}{| l | c | c | c |} 
         \hline
        & \textbf{Finite (unipotent)} & \textbf{Affine (universal)} & \textbf{Double affine (potent)} \\ 
         \hline
      
           \textbf{Weyl group} & $W$
            & $W\ltimes \Lambda$ &  $W\ltimes(\Lambda\times \check \Lambda)$ \\
           \hline

           \textbf{Hecke algebra} & $\HH^{\mathit{fin}}_G$ &  $\HH^{\mathit{aff}}_G$ &
           $\HH\HH_G$  \\
           \hline

           \textbf{Hecke ``category"} & $\cH_G^{\mathit{unip}} = \cD^{\mathit{unip}}(U\bs G/U)$ &  $\cH_G^{\mathit{univ}}=  \cD(U\bs G/U)$ &
           $\cH_G = {\mathbf D}^{pot}_{T\times T}(U\bs G/U) $  \\
           \hline
           

           \textbf{Representations} & $G-cat^{unip}$ &   $G-cat$ &
            $G-cat^{pot}$  \\
           \hline

           \textbf{Parameters} & $\check T^\wedge_e\GIT W$ & $\check T\GIT W$ & 
        $(\check T \times T)\GIT W$ \\
            \hline
    \end{tabular}

    \vspace{1em}
    
    \label{table}

\end{table}

\vspace{-2em}
The second row of traditional Hecke algebras  can be  concretely constructed  as deformations of the Coxeter data in the first row.  
The third row  contains monoidal Hecke categories  categorifyng the corresponding Hecke algebras in the case of the first two columns and conjecturally in the case of the third as discussed in~\S\ref{DAHA appears} below.
The fourth row contains the categorical (unipotent, de Rham or potent)
representation theory  governed by the corresponding Hecke category acting on the (unipotent, de Rham or potent) representation theory of the torus.
Finally, the fifth row lists the parameters of the Hecke category, or equivalently the representation theory it governs.  Note the dependence of $G$-categories 
on the parameters $\Tv\GIT W$ ~\cite{BZG} is completed/expanded in the unipotent/potent settings. We can also consider ``additive" versions, notably the graded potent theory $G-cat^{gr,pot}$ corresponding to the trigonometric Cherednik algebra, built on $W\ltimes (\Lambda \times \cO[\frt])$ and with parameters in $(\Tv \times \frt)\GIT W$.

Let us also highlight the effect of Langlands duality. On the one hand, in the first column, the finite Weyl group and Hecke algebra are Langlands self-dual, though at the categorical level we need to exchange monodromic and equivariant parameters. 
On the other hand, the third column is expected to be self-dual:
via the
decategorification discussed  in~\S\ref{DAHA appears} below,
the duality of the potent Hecke monad $\cH_G$ in our main Conjecture~\ref{potent duality conj}, a nonabelian refinement of the multiplicative 2-Fourier transform of Theorem~\ref{difference Fourier}, should lift the Langlands duality of the double affine Hecke algebra $\HH\HH_G$.

We should also make some clarifying comments: first, the potent Hecke monad $\cH_G = {\mathbf D}^{pot}_{T\times T}(U\bs G/U) $ is {\em not even a category}, but rather an object of an ind-coherent 2-category as constructed in \S\ref{GRT intro}. Consistent with this,
we do not know of a combinatorial construction of 
    $\cH_G $, though $\cH_G^{\mathit{unip}} $ and  $\cH_G^{\mathit{univ}}$ can be combinatorially constructed by Soergel bimodules.
Second, the universal Hecke category $\cH_G^{\mathit{univ}}$ is not the usual categorification of the affine Hecke algebra $\HH^{\mathit{aff}}_G$ given by the affine Hecke category $\cH_G^{\mathit{aff}} = \cD(I\bs LG/I)$, but provides a categorification viewing it through multiplicative Soergel bimodules~\cite{eberhardtK} or coherent Springer theory~\cite{BZCHN}. More generally we have minimized the role of loop groups and geometric Langlands  throughout. The tangled web of relations between the different realizations of these objects doesn't fit inside the physics of 4d gauge theory, but calls for 
the mysterious 6-dimensional (2,0) CFT (``Theory $\mathfrak X$''), where for example the category of representations of DAHA with all its symmetries appears via equivariant compactification on a four-torus(!).

\subsection{Taking traces: character sheaves, Hilbert schemes, DAHA, and link homology}\label{Hilbert section}
Here we discuss the decategorification of potent categorical representations -- so the evaluation of the 3d TQFT defined by $G-cat^{pot}$ on the circle $S^1$ -- and its expected relation to some active research directions. 

\subsubsection{Traces and character sheaves.}\label{character sheaves section}
Just as representations have characters, defined by taking traces of group elements, categorical representations, viewed as suitably dualizable objects of a 2-category $\cC$, have characters, which are objects of the trace (Hochschild homology) category $\Tr(\cC)$. 
Characters are invariant under the canonical $S^1$-action (cyclic symmetry) on the trace, and so further define objects of the periodic cyclic homology category $\Tr(\cC)^{\Tate}$, the Tate construction on the trace. 

When $\cC=A-\module$ for a monoidal category $A$, we have the trace category $\Tr(\cC)$ -- the ``horizontal decategorification" -- given by the cocenter of $A$, 
and the trace algebra $\tr(A)$ --
the ``vertical decategorification" -- given by taking
the trace (Hochschild homology) of $A$ as a plain category. In this case,  one finds $\tr(A)-\module$ as a full subcategory of $\Tr(\cC)$ by identifying $\tr(A)$ with endomorphisms of the trace of the monoidal unit of $A$, but $\Tr(\cC)$  is the richer object which we will focus on.

To approach the trace of $G-cat^{pot}$, let us review what is known about 
traces
for the first two columns of Table~\ref{table}.
In~\cite{character}, we calculated the cocenter 
(as well as the center, see also~\cite{BFO})
of the unipotent Hecke category $\cH_G^{\mathit{unip}}$, and 
thus the trace of unipotent categorical representations
$G-cat^{\mathit{unip}}$,
to be Lusztig's~\cite{lusztigcharacter}  {\em unipotent character sheaves} $\Ch^{\mathit{unip}}_G \subset\cD(G/G)$,  class $\cD$-modules with nilpotent singular support and unipotent central character.
In particular, 
thanks to the Koszul duality of~\cite{BGS}, this identifies  the categories of unipotent character sheaves for $G$ and $\Gv$, 
an enrichment of the identification of the correpsonding finite Hecke algebras here encoded by Springer blocks. 
In~\cite{characterfieldtheory}, we showed that passing from  $G-cat^{\mathit{unip}}$ to 
all de Rham categorical representations $G-cat$ 
leads on the level of traces to all class $\cD$-modules $\cD(G/G)$.

\subsubsection{Context: Link homology and Hilbert schemes} 
Before turning to the trace of $G-cat^{pot}$, we provide some further geometric context and motivation here.
A major impetus for 
Hecke categories (or Soergel bimodules) and their traces comes from 
Khovanov-Rozansky link homology theory~\cite{KhovanovRozansky,Khovanovtriplygraded}. The Artin braid group of $G$ sits inside the finite Hecke category, whence applying various functors gives homological invariants of braids. 
In particular, the 
link homology invariant of a braid can be obtained~\cite{WebsterWilliamson}  by taking the cohomology of the 
trace of the corresponding object in the (mixed) Hecke category for $GL_n$. In this way, the link homology construction is a functor of the trace of the braid regarded as a unipotent character sheaf~\cite{HoLiHilbert}. 

On the other hand, a remarkable circle of ideas has emerged relating Khovanov-Rozansky link homology to coherent sheaves on the Hilbert scheme  $\textup{Hilb}_n(\CC^2)$ of points on the plane and representations of rational Cherednik algebras ~\cite{GORS,GNR,OblomkovRozanskyHilbert}. However, the sheaves on $\textup{Hilb}_n(\CC^2)$ appearing in these conjectures are supported on the locus ``$y=0$'', i.e.~the subscheme of points supported on the $x$-axis in the $x,y$-plane. Indeed, a theorem of~\cite{HoLiHilbert} identifies the trace of the mixed finite Hecke category
with equivariant quasicoherent sheaves on this subscheme. The ``$y$-ification'' program begun in~\cite{GorskyHogancamp} proposes a deformation of the Hecke category, and thereby of the associated link homology theories,  off of the $x$-axis to a theory
which is Koszul self-dual (exchanging equivariant and monodromic -- or $x$ and $y$ -- parameters) and whose trace maps to (but differs from) quasicoherent sheaves on the full Hilbert scheme. The self-duality is known to be related to fundamental symmetry in combinatorics, representation theory of Cherednik algebras and Hodge theory of character varieties~\cite{GorskyHogancampMellit}.

\subsubsection{Traces of potent representations}\label{DAHA appears}
We suggest $y$-ification is profitably viewed through the lens of potent representation theory. We expect the trace of potent representations will see the entire Hilbert scheme, while also providing multiplicative analogues.

To start, in the case of a torus $T$,  we have discussed 
that 
$T-cat^{pot, gr}$ 
models the Rozansky-Witten theory with target  $\frt\times \check T$ 
with trace $\QCoh(\frt\times \check T)$. Thus for the maximal torus $T \subset GL_n$, the trace lives on the $n$-fold product of  $\CC \times \Cx$, 
and in particular the $y$-parameters appear for free. 
From this starting point, following the results of~\cite{HoLiHilbert} (and the McKay correspondence) for unipotent Hecke categories, we are led to the expectation that the trace of $GL_n-cat^{pot, gr}$ 
is a form of $\QCoh(\textup{Hilb}_n(\CC\times \Cx))$
(the resolved Coulomb branch).

In the (non-exact!) multiplicative setting of $T\times \Tv$, the definition of extended RW theory (even where it is valued) is far more subtle, but we sketch some expectations. We propose $T-cat^{pot}$ as a model for RW theory with target $T\times \Tv$ with  trace a form of  $\QCoh(T\times \Tv)$. 
For $G=GL_n$, we  expect the trace of $GL_n-cat^{pot}$ to be a form of $\QCoh(\textup{Hilb}_n((\Cx)^2)$ (the resolved K-theoretic Coulomb branch) compatibly with Langlands self-duality. 


What about the cyclic trace $Tr(G-cat^{pot})^{Tate}$ (or better, from the 4d point of view, $\cB_G(T^2)^{T^2-\Tate}$)? 
As we discussed in~\S\ref{Omega section}, the $\Omega$-background implements deformation quantization of the moduli of vacua. In the nonexact setting, this depends on a parameter $q$ which for us arises from the ratio of two $S^1$-equivariant parameters. For a torus $T$, this results in the appearance of categories of difference modules on $T$. For a general reductive group $G$, we expect to see DAHA
specialized at $t=1$.
More specifically, the theory of quantum groups provides the category $\cD_q(G/G)$  of equivariant quantum $\cD$-modules~\cite{BZBJ2,BrochierJordan,quantumcharactertheory}. It contains a Springer block which is equivalent (in type $A$) with modules for the $t=1$ specialization of DAHA, 
as well as a ``difference Fourier transform'' autoequivalence restricting to that of DAHA. For $G=GL_n$, we expect a relation between cyclic traces of potent $G$-categories and $\cD_q(G/G)$ compatible with the Fourier transform, i.e.~with representations of DAHA.  Thus potent $GL_n$-categories should define DAHA modules as characters, compatibly with Langlands/Fourier duality. We further hope a mixed version of potent representations in the spirit of~\cite{HoLiMixed} sees the full $(q,t)$ DAHA theory. 



\subsection{The character TQFT and potent homology of character stacks}\label{surfaces section}
Here we discuss the verification that  $G-cat^{pot}$ defines a (partial) 3d TQFT, and its  evaluation
on surfaces, specifically in  relation to the homology of character stacks.


 
 In joint work with  Gunningham~\cite{character},~\cite{characterfieldtheory} we applied the Cobordism Hypothesis~\cite{jacobTFT} to construct  [0,2]-extended oriented 3d TQFTs $\Xi^{\mathit{unip}}_G$ and $\Xi_G$ with the assignments
 $\Xi^{\mathit{unip}}_G(pt) = G-cat^{\mathit{unip}}$ and $\Xi_G= G-cat$. We use the name {\em unipotent character theory} and
 {\em character theory} for $\Xi^{\mathit{unip}}_G$ and $\Xi_G$ due to their assignment of the respective character sheaves to the circle $S^1$ as recalled in~\S\ref{character sheaves section} above.
For a closed oriented surface $C$, $\Xi_G(C)$ produces the (renormalized Borel-Moore) homology of the character stack $\Loc_G(C)$. On the other hand by Koszul duality~\cite{BezYun} the unipotent homology  $\Xi_G^{\mathit{unip}}(C)$ is identified for Langlands dual groups.

We conjecture that the richer 2-category $G\cat^{pot}$ is likewise 2-dualizable and 2-orientable, so that we may plug it into the Cobordism Hypothesis and build a (partial) 3d TQFT $\Xi_G^{pot}$ called the {\em potent character theory}. It would provide a model for Witten's equivariant compactification of Kapustin-Witten theory as discussed in ~\S\ref{4d section}.
The duality of Conjecture~\ref{potent duality conj} would immediately imply a Langlands duality for the value of this TQFT on oriented surfaces:
$\Xi_G^{pot}(C)\simeq \Xi_{\Gv}^{pot}(C)$.
We refer to this assignment as the {\em potent homology} of character stacks 
and expect there is a description of it in terms of classical invariants of character stacks such as equivariant K-theory. 

Finally, following ~\S\ref{Hilbert section}, for $G=GL_n$, recall we expect 
the value $\Xi_G^{pot}(S^1)$ on the circle to give quasicoherent sheaves on the Hilbert scheme of $(\Cx)^2$. It follows that after fixing a point $c\in C$,  the potent homology of the character stack $\Xi_G^{pot}(C)$ would be the global sections of a quasicoherent sheaf $\Xi_G^{pot}(C\setminus x)$ on the Hilbert scheme. This would be compatible with Langlands duality, and we hope shed light on the deep relations between homology of character varieties, Macdonald polynomials, and DAHA~\cite{HRV, MellitECM}. 


\subsection{Potent relative, real and finite Langlands correspondences}
We conclude here by discussing important examples of objects in 
$G-cat^{pot}$ and their predicted Langlands duality. 
Recall from Example~\ref{constructing potent reps} that  a $G$-space $X$ defines a potent $G$-category 
${\mathbf D}^{pot}_G(X)$, which then has a character
in the trace of $G-cat^{pot}$ (so for $G=GL_n$ conjecturally a quasicoherent sheaf on the Hilbert scheme of $(\Cx)^2$) and in the cyclic trace (whence conjecturally a module for DAHA).

\subsubsection{Relative Langlands}
The relative Langlands duality developed in~\cite{BZSV} provides a rich source of matching measurements on the two sides of the Langlands correspondence, extending the role of L-functions in the classical theory of automorphic forms. It is based on the symplectic and symmetric perspective of TQFT, specifically the S-duality of boundary theories in 4d gauge theory~\cite{GaiottoWittenSduality}. Applying equivariant compactification on a circle, relative Langlands suggests many examples of dual pairs of boundary theories for the potent character TQFT, i.e.~dual pairs of potent categorical representations.

Let $G\actson M$ be a hyperspherical $G$-variety and $\Gv\actson \Mv$ its Langlands dual hyperspherical $\Gv$-variety. We assume for simplicity that $M$ and $\Mv$ are both polarized, $M=T^*X$ and $\Mv=T^*\Xv$ for suitable $G\actson X$ and $\Gv\actson \Xv$ (one can also twist to include Whittaker type examples). The key compatibility between these dual spaces is the local unramified geometric conjecture of \cite{BZSV} which identifies sheaves on $LX/LG_+$ with $\IndCoh(\cL \Xv/\Gv)$ as modules for the spherical Hecke category. This conjecture (which is known in many cases) has natural tamely ramified extensions to module categories for the affine Hecke category (see e.g.~\cite[Conjecture 3.4.14]{Devalapurkar}). Applying the Tate construction for loop rotation is expected~\cite[Remark 1.1.4]{FGT} to lead to a Koszul duality \cite[Conjecture 1.1.3]{FGT} between the categories of equivariant $\cD$-modules $\cD(X/B)$ and $\cD(\Xv/\Bv)$ on dual spherical varieties as modules for the finite Hecke category. This unipotent conjecture has a natural potent upgrade: 

\begin{conj}[Potent relative Langlands duality]\label{potent relative}
The potent Langlands duality of Conjecture~\ref{potent duality conj} matches ${\mathbf D}^{pot}_G(X)$ with 
${\mathbf D}^{pot}_\Gv(\Xv)$ for dual polarized hyperspherical varieties.
\end{conj}



\subsubsection{Real local Langlands}
One of our original motivations for many of the constructions discussed is the real local Langlands correspondence, in the form developed by~\cite{ABV,Soergel} and reformulated in~\cite{reps,BZCHN2} via the geometry of derived loop spaces.

 Let $\theta$ be a quasisplit real form of $G$, and $\check \theta$ the dual involution of $\Gv$. They give rise to collections $\Theta$ of pure inner forms of $G$ and $\check\Theta$ of involutions of $\Gv$, which label the stacky fixed points on the level of classifying spaces, $BG^{\theta}=\coprod_{\Theta} BG_\tau$ and likewise for $\check\Theta$.
Soergel conjectured (with the quasi-split case proved~\cite{BezVil}) a Koszul duality between the categories of Harish-Chandra modules for all $G_\tau$ with trivial infinitesimal character and sheaves $\bigoplus \cD(\Bv\bs \Gv/\Gv_{\check \tau})$, compatibly with actions of the finite Hecke category, recovering a theorem of~\cite{ABV} on the level of Grothendieck groups. More generally as we vary the infinitesimal character for $G$ we replace $\Gv$ by a corresponding semisimple centralizer. 

In~\cite{reps,BZCHN2}) this family of conjectures was reorganized using Jordan decomposition of loop spaces and potent $\cD$-modules as in \S\ref{loops and conns section}. Namely, we showed that the Soergel conjecture for arbitrary regular infinitesimal character can be deduced (by completing at a parameter) from a uniform statement ~\cite[Conj.4.5]{BZCHN2}.
This conjecture in turn would follow via the Tate construction from a geometric Langlands conjecture on the twistor $\PP^1$, in agreement with Scholze's recent geometrization of real local Langlands~\cite{scholzetwistor}.

We now suggest a stronger and more symmetric potent form of real local Langlands: 

\begin{conj}[Potent Soergel-Vogan duality]\label{potent Soergel}
The potent Langlands duality of Conjecture~\ref{potent duality conj} matches 
 the potent $G$-category represented
by loops on $(BG)^\theta$ with the potent $\Gv$-category represented by loops on~$(B\Gv)^{\check{\theta}}$.
\end{conj}

Passing to characters, we thus expect a {\em potent Vogan duality}, a duality of representations of double affine Hecke algebras expanding Vogan's duality of K-groups of representations of a fixed infinitesimal character as modules for the finite Hecke algebra. 


\subsubsection{Potent representations of finite reductive groups}
While we have worked throughout over the complex numbers, it is expected that many of the constructions have analogues in positive characteristic. The modern theory of categorical traces of Frobenius~\cite{dennisshtuka,AGKRRV1, xinwentrace} was applied in~\cite{etevetrace} to recover and generalize aspects of Lusztig's theory of representations of finite reductive groups from Koszul duality for finite Hecke categories. By applying the trace of Frobenius to endoscopic Koszul duality results~\cite{BezYun,LusztigYun}, one obtains equivalences of categories between representations of $G(\FF_q)$ with semisimple parameter $\sv\in \Gv(\FF_q)$ and unipotent representations of the centralizer of $\sv$, or equivalently of its Langlands dual group, the endoscopic group of $G(\FF_q)$ associated to $\sv$. 

We expect there is a potent form of categorical representations of $G/\overline{\FF_q}$. One approach would be to apply the Tate construction (for rotations of the punctured disc) to ind-coherent sheaves of categories on the stack of restricted~\cite{AGKRRV1} tamely ramified local Langlands parameters, cf. Gaitsgory's article in these proceedings.  One could then define the category of potent representations of $G(\FF_q)$ as the trace of Frobenius on potent categorical representations over $\overline{\FF_q}$, and conjecture that the categories of potent representations of Langlands dual Chevalley groups $G(\FF_q)$ and $\Gv(\FF_q)$ are equivalent.
Evidence for this is available on the semisimple level of complex representations, where it amounts to the following observation: Lusztig's classification becomes Langlands self-dual if we replace representations of $G(\FF_q)$, which depend on a semisimple parameter $\sv\in \Tv$, by the sum over pairs $s,\sv\in T\times \Tv$ (equipped with $\FF_q$-rational structure) 
of {\em unipotent} representations of the ``joint centralizer'' $G_{s,\sv}(\FF_q)$, the centralizer of $s$ in the endoscopic group for $\sv$, or equivalently, unipotent representations of the dual group $\Gv_{s,\sv}$.





\bibliographystyle{alphaurl}
\bibliography{ICM}

\end{document}